\documentclass[a4paper, 10pt]{article}

\usepackage[a4paper,left=3cm,right=3cm,top=2.5cm,bottom=2.5cm]{geometry}
\usepackage{hyperref}
\usepackage{subfig}
\usepackage{pgfplots}
\usepackage{pgfplotstable}
\usepackage{mathtools}
\usepackage{multicol}
\usepackage{comment}
\usepackage{booktabs}
\pgfplotsset{compat=1.5}
\usepackage{amssymb}
\usepackage{url}
\usepackage{bm}

\usepackage{pdfsync}
\usepackage{float}
\usepackage{tabularx}
\usepackage{enumerate}
\usepackage{array}
\usepackage{xspace}
\usepackage{tikz}
\usepackage{tikz-cd}
\usepackage{tikzsymbols}
\usetikzlibrary{calc,trees,positioning,arrows,chains,shapes.geometric,%
    decorations.pathreplacing,decorations.pathmorphing,shapes,%
    matrix,shapes.symbols, decorations.markings, patterns,fit}

\usepackage{siunitx}

\usepackage{authblk}

\usepackage[draft,inline,marginclue]{fixme}
\usepackage{mathrsfs}
\usepackage{color, colortbl}
\usepackage{multirow}

\usepackage{amsthm}
\usepackage{amsmath}
\usepackage{stmaryrd}
\usepackage{algorithm}
\usepackage{algpseudocode}
\usepackage{graphicx}
\usepackage{booktabs}

\usepackage[noabbrev]{cleveref}

\newcommand{\mupar}{\ensuremath{\boldsymbol{\mu}}}
\newcommand{\etapar}{\ensuremath{\boldsymbol{\eta}}}
\newcommand{\x}{\ensuremath{\mathbf{x}}}
\newcommand{\s}{\ensuremath{\mathbf{s}}}
\newcommand{\co}{\ensuremath{\mathbf{c}}}

\newcommand{\R}{\ensuremath{\mathbb{R}}}
\DeclareMathOperator*{\argmin}{\arg\!\min}
\DeclareMathOperator*{\argmax}{\arg\!\max}
\newcolumntype{C}[1]{>{\centering\arraybackslash}m{#1}}

\definecolor{Gray}{gray}{0.9}

\newcommand{\RA}[1]{{\color{black}#1}}
\newcommand{\RB}[1]{{\color{black}#1}}

\begin{document}

\title{A multi-fidelity approach coupling parameter space reduction
  and non-intrusive POD with application to structural optimization of passenger ship hulls}

\author[a]{Marco~Tezzele\footnote{marco.tezzele@sissa.it}}
\author[a]{Lorenzo~Fabris\footnote{lorenzo.fabris@sissa.it}}
\author[b]{Matteo~Sidari\footnote{matteo.sidari@fincantieri.it}}
\author[b]{Mauro~Sicchiero\footnote{mauro.sicchiero@fincantieri.it}}
\author[a]{Gianluigi~Rozza\footnote{gianluigi.rozza@sissa.it}}

\affil[a]{Mathematics Area, mathLab, SISSA, via Bonomea 265, I-34136 Trieste,
  Italy}
\affil[b]{Fincantieri S.p.A., Merchant Ships Business Unit, Passeggio
  Sant'Andrea 6/A, I-34123 Trieste, Italy} 

\maketitle

\begin{abstract}
Nowadays, the shipbuilding industry is facing a radical change towards
solutions with a smaller environmental impact. This can be achieved
with low emissions engines, optimized shape designs with lower wave
resistance and noise generation, and by reducing the metal raw
materials used during the manufacturing. This work focuses on the last
aspect by presenting a complete structural optimization pipeline for modern
passenger ship hulls which exploits advanced model order reduction
techniques to reduce the dimensionality of both input parameters and
outputs of interest. We introduce a novel approach which incorporates
parameter space reduction through active subspaces into the proper
orthogonal decomposition with interpolation method. This is done in a
multi-fidelity setting. We test the whole framework on a simplified model of a midship section and 
on the full model of a passenger ship, controlled by $20$ and $16$ parameters, 
respectively. We present a comprehensive error analysis and show the
capabilities and usefulness of the methods especially during the
preliminary design phase, finding new unconsidered designs while
handling high dimensional parameterizations.
\end{abstract}

\tableofcontents

\section{Introduction}
\label{sec:ironth_intro}
When considering optimization of complex systems in an industrial
context we must rely on surrogate models in order to alleviate the
computational cost of this kind of many-query
problems~\cite{chinestaenc2017, beohparour17}.
Scientific machine learning~\cite{baker2019workshop} is widely used in
applied mathematics and in engineering
applications~\cite{ghattas2021learning, brunton2021data, brunton2019data} such as inverse problems, optimization, and
prediction of the behaviour of parametrized systems, to cite a few. In
this work we are considering structural reduced order models (ROMs) of
passenger ship hulls in order to speed up the optimization process.
Structural analysis of complex systems through reduced order modeling
is not limited to naval engineering. Recently a component-based
data-driven approach has been proposed to assess the structural integrity of aircraft
components~\cite{kapteyn2020data, kapteyn2021probabilistic} in the context of modern digital
twins incorporating not only data but also physical models, also
referred to as hybrid twins~\cite{chinesta2020virtual}.
For multidisciplinary analysis and optimization involving reduction in
both input and output spaces we cite~\cite{boncoraglio2020model,
  boncoraglio2021active}, while for a specific naval engineering
application we suggest~\cite{khan2021physics}.

In this work we propose an optimization framework, to be used in the
preliminary design phase, involving
many reduced order models to assess the structural behaviour of modern
passenger ship hulls under different parametric configurations and
loading conditions.
Many studies have been conducted to assess the structural behaviour of
passenger ship hulls~\cite{heder1991hull, soares2017progress,
  romanoff2013hull, parmasto2012mechanics}. 
In~\cite{prebeg2014application} they compare different surrogate
models to improve the design process of complex thin-walled ship
structures, without using any \RA{proper orthogonal decomposition (POD)}-based model order reduction.
For structural behaviour and optimization of passenger ships we
cite~\cite{raikunen2019optimisation},
and~\cite{raikunen2018optimization} where they used
efficient finite element modelling, evolutionary optimization
algorithm and indirect constraint relaxation. We used a similar idea
for the stability constraints, where local stress peaks
are allowed to exceed the rule-based strength limits.

\RA{A recent approach to accelerate PDE-constrained optimization was proposed
in~\cite{keil2022adaptive}, where an adaptive method comprising both
full order model evaluations and artificial neural networks surrogate
models was used in the context of oil recovery. The main idea
introduced was to use local approximations of the objective functional
instead of a global surrogate model.
Another adaptive numerical method involving ROMs and nonlinear
trust-region based on a residual error indicator able to keep the
optimization trajectory consistent with the ROMs accuracy was presented
in~\cite{zahr2015progressive}.  
}

The novelty of this work is the incorporation into the proper
orthogonal decomposition framework of parameter space reduction~\cite{tezzele2022aroma16} by
constructing a multi-fidelity surrogate model~\cite{romor2021multi, perdikaris2017nonlinear}, without the need of
running simplified simulations. This is done by exploiting the presence of an active
subspace~\cite{constantine2015active} of the parameter to reduced
state variables map. This represents
a new data-driven non-intrusive ROM, more accurate with respect
to a more classical interpolation method such as Gaussian process
regression (GPR)~\cite{williams2006gaussian}. \RA{We called this new
  method POD-NARGPAS since it comprises POD with interpolation,
  nonlinear autoregressive Gaussian processes, and active
  subspaces. We also introduce a structural optimization numerical pipeline for
  large scale industrial applications, exploiting the newly proposed
  ROM.
We use a Bayesian approach to perform discrete mono-objective
structural optimization. We incorporate stability constraints in a
weak form by accounting for the added mass needed to stabilize the
affected elements. 
The optimization pipeline, thanks to its modularity, allows for different
target functions to minimize, from the total mass of the hull, to the manufacturing
cost of the structure. \RA{Moreover, it has the potential of being
  used in many other engineering fields for surrogate-based optimization
  tasks~\cite{forrester2009recent, koziel2011surrogate}, especially
  the ones involving costly numerical simulations.} 

In naval engineering multi-fidelity methods have been used in the
context of surrogate-based design optimization for super-cavitating
  hydrofoils~\cite{bonfiglio2018multi}, marine propellers
  design~\cite{gaggero2022marine}, and the optimization of a NACA
  hydrofoil with an adaptive sampling method using stochastic radial
  basis functions~\cite{pellegrini2022towards}, for example.}

This work is organized as follows: in
\autoref{sec:ironth_opt_pipeline} we present the entire pipeline; in
\autoref{sec:ironth_fom} and \autoref{sec:ironth_rom} we describe the full
order model and the reduced order ones, respectively. We introduce
parameter space reduction with active subspace, proper orthogonal
decomposition with interpolation, and how to combine them in a
multi-fidelity autoregressive scheme. In \autoref{sec:ironth_bayes} we briefly summarize the
Bayesian optimization scheme we used for the numerical results reported in
\autoref{sec:ironth_results}, where we tested the framework on a midship
section of a simplified hull and on a real complete hull. Finally in
\autoref{sec:ironth_conclusions} we draw the conclusions and some future
research lines.

\section{Structural optimization pipeline}
\label{sec:ironth_opt_pipeline}

In this work we use the Nested Analysis and Design
(NAND)~\cite{arora2005review, balesdent2012survey} approach. We
consider the problem:
\begin{equation}
\min_{\substack{\mupar \in \mathcal{P}}}  \quad f(\s, \mupar),
\qquad \text{s.t.} \quad R(\s, \mupar) = 0, 
\end{equation}
where $\s$ represents the state vector, and $R(\s, \mupar)$ a general
high-dimensional discretized parametric \RA{partial differential equation
(PDE)}, which we are going to characterize in
\autoref{sec:ironth_fom}. 
Following the NAND approach, where $\s$ is considered an
implicit function of $\mupar$, we can rewrite the optimization problem as
\begin{equation}
\min_{\substack{\mupar \in \mathcal{P}}}  \quad f(\s (\mupar), \mupar).
\end{equation}
So for every queried parameter point $\mupar$, we solve the PDE and
evaluate the function to minimize. For a fast and accurate solution of
the PDE we use reduced order models described in \autoref{sec:ironth_rom}.

\begin{figure}[htb]
\centering
\includegraphics[width=1.\textwidth]{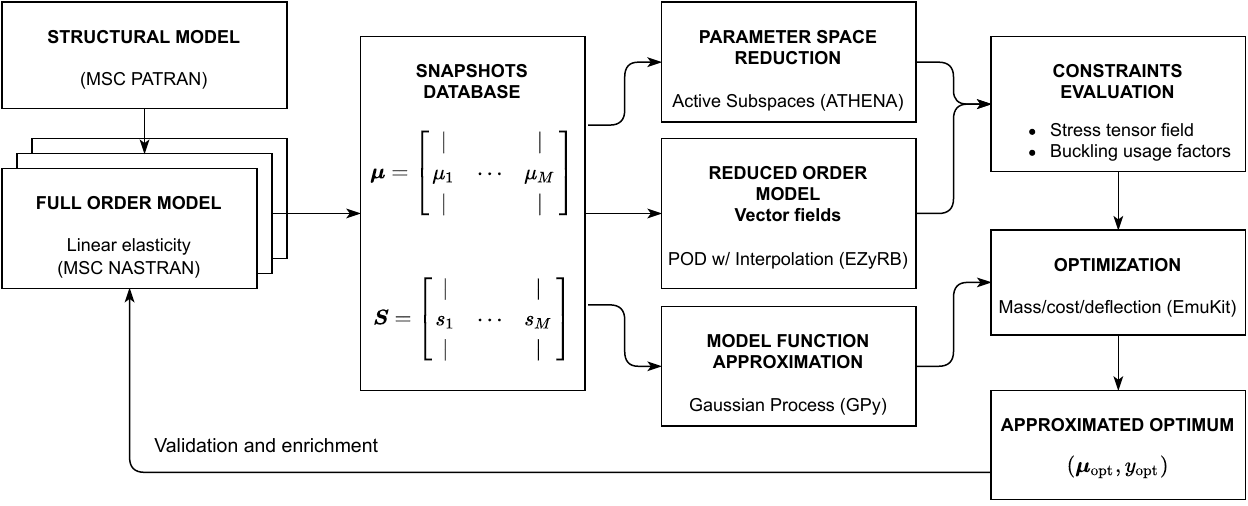}
\caption{Structural optimization workflow, from the base structural model
  creation, to the approximated optimum and validation. Each block of
  the pipeline reports the underlying method and the software
  used.}
\label{fig:ironth_hull_opt_scheme}
\end{figure}

The complete structural mono-objective optimization workflow is depicted in
\autoref{fig:ironth_hull_opt_scheme}, in which for every building block we
emphasize the software used.

We start from the construction of the parametrized structural model
with MSC Patran, and we construct a database of full order solutions
with MSC Nastran corresponding to a given set of parameters for every
loading condition. With this database we construct different 
reduced order models depending on the quantity of interest we
want to approximate. We use POD with interpolation (PODI)~\cite{ly2001, bui2003proper} for
the stress tensor field approximation, and
GPR~\cite{williams2006gaussian} for the approximation of scalar
functions.
Moreover we exploit \RB{active subspaces~\cite{constantine2015active} (AS)} for the reduction of
the parameter space dimension to build low-fidelity models and improve
the PODI prediction capabilities in a multi-fidelity
setting~\cite{peherstorfer2018survey, perdikaris2017nonlinear} called nonlinear
autoregressive multi-fidelity Gaussian process regression with active
subspaces (NARGPAS)~\cite{romor2021multi}. These
parameter and model reduction methods are combined for a computational
efficient and reliable evaluation of the constraints regarding the
stability of the whole hull: we check how many elements are yielded,
and how many elements are subjected to buckling phenomena. We remark
that we allow local stress peaks to exceed the classification society rule
limits, since we automatically incorporate within the function
to optimize the necessary interventions at the shipyard to stabilise
such elements. This is particularly important since the proposed
pipeline is going to be used in the preliminary design phase. The
optimization is done with a Bayesian approach. The approximated
optimum is then validated with the full order model and the snapshots
database enriched accordingly.

\begin{figure}[htb]
\centering 
\includegraphics[trim=0 0 20 60, clip, width=.49\textwidth]{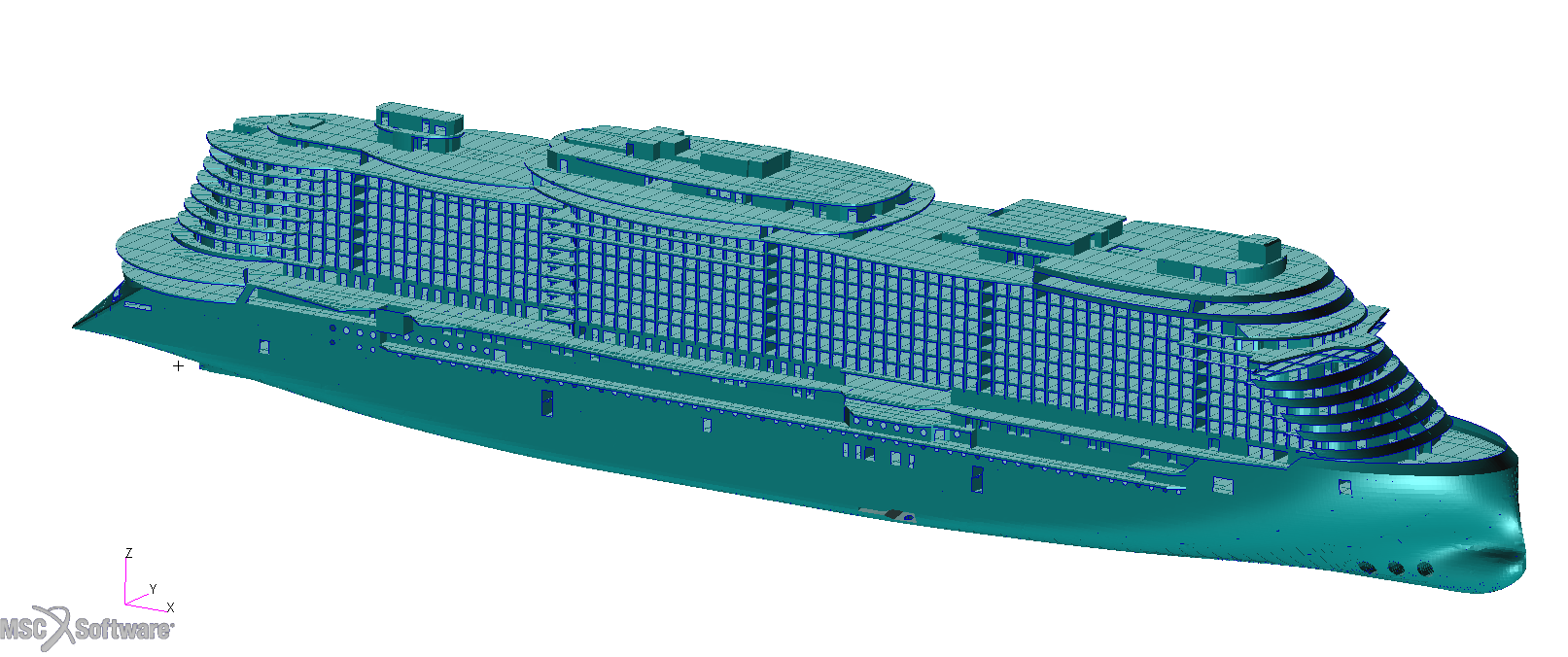}\hfill
\includegraphics[trim=0 0 20 60, clip, width=.49\textwidth]{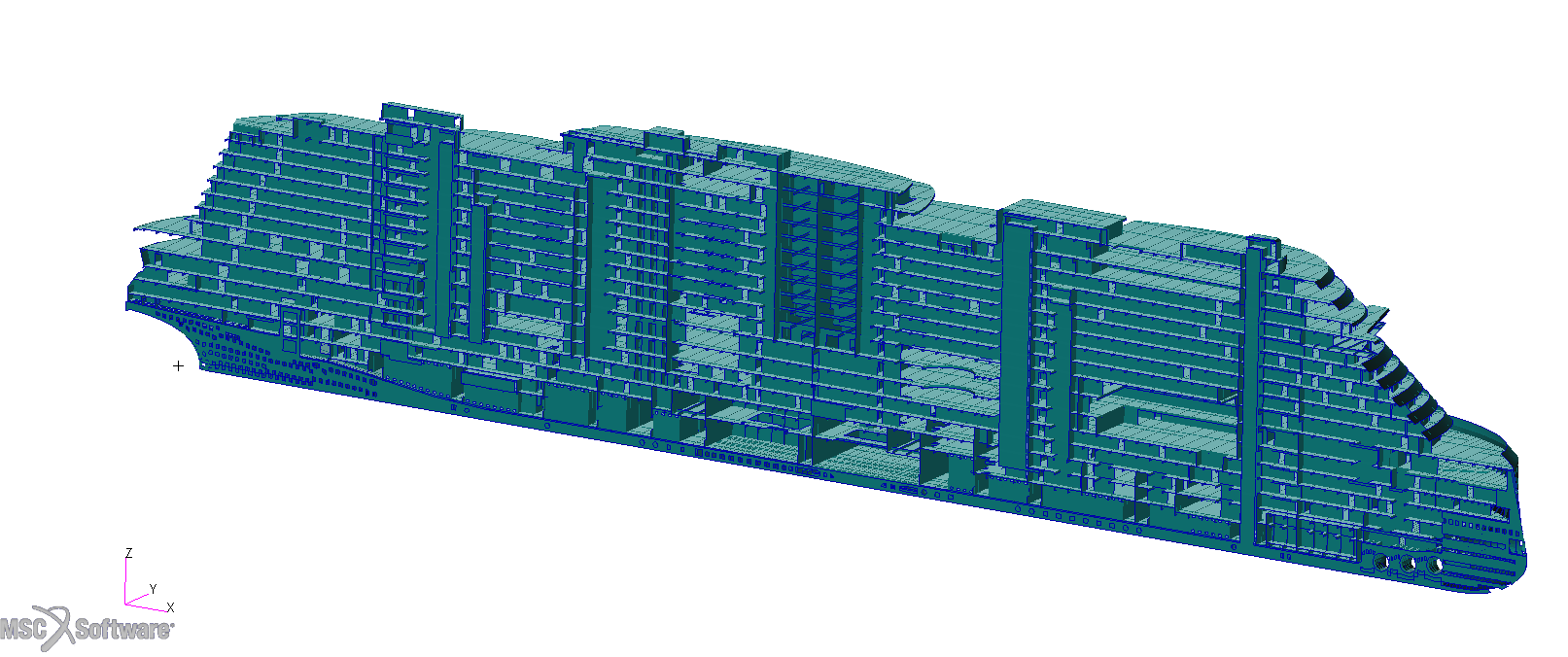}
\caption{A complete view of the hull on the left, and a longitudinal section on the right.}
\label{fig:ironth_c6298}
\end{figure}

We are going to present all the numerical methods employed and finally
the application of the whole pipeline to an actual passenger ship
hull depicted in \autoref{fig:ironth_c6298} and built by Fincantieri S.p.A.. 

\section{Full order model}
\label{sec:ironth_fom}

In this section we describe the PDE we need to solve and the
high-fidelity solver used to create the solutions database.

The equations governing the linear elastic isotropic problem are the
equilibrium equation, the linearised small-displacement
strain-displacement relationship, and the Hooke's law, respectively:
\begin{equation}
  \begin{cases}
 & - \nabla \cdot \sigma = h,\\
 & \epsilon = \frac{1}{2} [ \nabla u + \nabla u^T ], \\
 & \sigma = C(E, \nu) : \epsilon,
\end{cases}
\end{equation}
where $\sigma$ is the Cauchy stress tensor, $h$ is the body force,
$\epsilon$ is the infinitesimal strain tensor, $u$ is the displacement
vector, and $C$ is the fourth-order stiffness tensor depending on $E$,
the Young modulus, and on $\nu$, the Poisson's ratio.

\begin{figure}[h!]
\centering
\includegraphics[width=.6\textwidth, trim=0 0 95 105, clip]{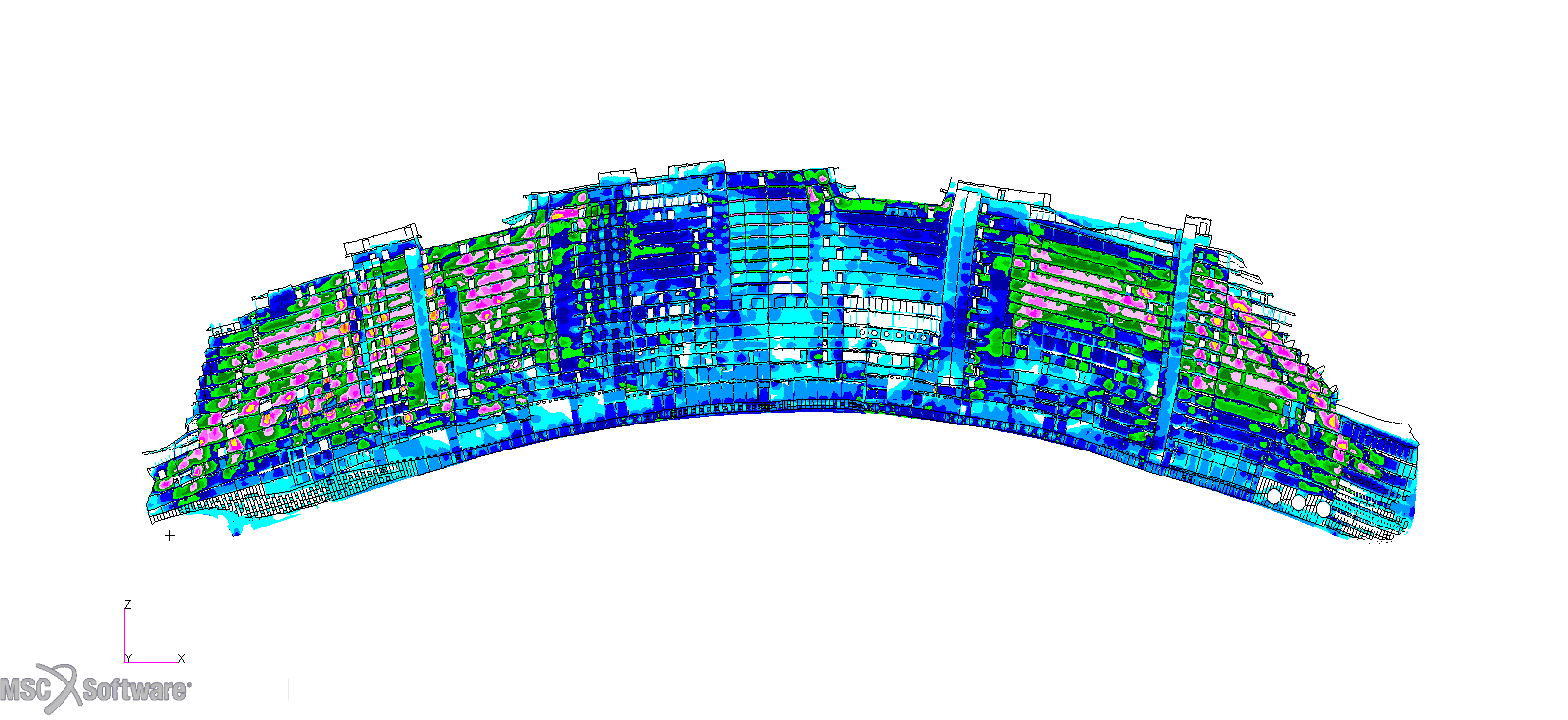}
\caption{Possible deformation of the hull under the hogging
  loading condition. Displacements are magnified. Colors refer to the
  von Mises criterion.}
\label{fig:ironth_c6298_vm}
\end{figure}

The finite element method employed uses $2$-dimensional elements,
commonly referred to as plate and shell elements. They are
used to represent areas in the model where one of the dimensions is small in
comparison to the other two. The
height or thickness of the element is substantially less than the
width and the length. We are going to use MSC Nastran CQUAD$4$ and
CTRIA$3$ elements, which are general-purpose plate elements capable of carrying
inplane force, bending forces, and transverse shear force. The
membrane stiffness of the $2$-dimensional elements is calculated using
the plane stress theory. Most thin structures constructed from common
engineering material, such as aluminum and steel, can be modeled
effectively using plane stress. In this work we consider
only high strength structural steel (AH26) for illustrative
reasons. The parameterization consists of the
thickness associated to specific regions of the hull. We compute the
solution for two classical loading conditions, namely the hogging and
sagging. In \autoref{fig:ironth_c6298_vm} an
example of solution for the hogging loading condition. We are
interested in the stress tensor field, through which we can compute
the von Mises criterion and the buckling usage factors, used for the
constraint's evaluation.

\section{Reduced order models}
\label{sec:ironth_rom}

In this section we describe the new proposed non-intrusive
data-driven ROM exploiting the low-intrinsic dimensionality of
the parameters to reduced snapshots map.
In this work the parameters vector $\mupar \in \mathcal{P}$
represents the thickness of some selected macro areas of steel plates.
In order to speed up the optimization procedure we construct a reduced
order model for all the stress tensor components, and from them we compute the
derived quantities which describe the constraints and the functions to
optimize.

\subsection{Sampling strategy}
\label{sec:podi}
Since the steel plates can have only a finite set of possible
thickness, the input parameter space is discrete. In
order to cover the domain in a uniform unbiased way, we generate many
random uniform sampling $\mathcal{P}^i := \{ \mupar^i_j\}_{j=1}^M$ composed by $M$
samples each. We compute the minimum pairwise distance $d$
associated to each $\mathcal{P}^i$ and we retain the samples with maximum
$d$. Let us call $\mathcal{P}$ the chosen samples set $\mathcal{P} := \{
\mupar_j\}_{j=1}^M$. \RB{This sampling plan allows to achieve univariate uniform distribution
without manipulation of a Latin hypercube, while maximizing a space-filling metric
akin to the maximin Latin hypercube design technique
from~\cite{johnson1990LHD}. More sophisticated sampling schemes for
discrete coordinates can be found in~\cite{you2021mp}.}

\subsection{Proper orthogonal decomposition with interpolation}

Non-intrusive data-driven POD-based reduced order models are reviewed
in~\cite{tezzele2022aroma09} with different applications. Here we
briefly present POD with interpolation.

For every $\mupar \in \mathcal{P}$ we solve the associated high-fidelity
parametric problem defined above, and we store the solution snapshots
$\s_j := \s (\mupar_j) \in \R^n$, with $j \in [1, \dots, M]$, in matrix form as
follows:
\begin{equation}
  S=\left[
    \begin{array}{cccc}
      | & | &  & |\\
      \s_1 & \s_2 & \dots & \s_M\\
      | & | &  & |
    \end{array}
  \right].
\end{equation}
We assume the state variables can be approximated as a linear combination
of a few global basis functions, also called modes, that is
\begin{equation}
\s_i = \sum_{k=1}^M \varphi_k  c_i^k \approx \sum_{k=1}^r \varphi_k  c_i^k,
\qquad \forall i \in [1, \dots, M], \;\; r \ll M,
\end{equation}
for some modes $\varphi_k \in \mathbb{R}^n$ and for some modal
coefficients $c_i^k \in \mathbb{R}$. Equivalently in matrix form we
have $S \approx \Phi C$, with $S \in \mathbb{R}^{n \times
  M}$, $\Phi \in \mathbb{R}^{n \times r}$, and $C \in \mathbb{R}^{r
  \times M}$. To compute such modes we use the
proper orthogonal decomposition technique. We decompose the
matrix $S$ with the truncated singular value decomposition:
\begin{equation}
S = U \Sigma V^* \approx U_r \Sigma_r V_r^*, \;\; r \ll M,
\end{equation}
where $^*$ denotes the conjugate transpose and the
subscript $r$ indicates the first $r$ columns. The columns of $U$
span the optimal low-dimensional subspace in the least square sense
and are also called POD modes. We have $\Phi := U_r$. 
To compute the modal coefficients $C$, also called reduced state
variables, we project the data onto the POD subspace: $C = \Phi^T S$.

With the matrices $\Phi$ and $C$ we are able to reconstruct the
initial database of solutions in a reduced way, but we are not able to
predict the state vector corresponding to a new parameter
$\mupar^*$. In order to do so, we need to compute the parameters to
reduced states map and then use the global basis to predict the
entire stress field. We construct a function $g : \mathcal{P} \to \R^r$ which approximates the
map $\mupar \in \mathcal{P} \to \co \in \R^r$, given the initial set of $M$ input-output
pairs $\{\mupar_j, \co_j\}_{j=1}^M$, where $\co_j = \Phi^T \s (\mupar_j)$. The
regression model $g$ is used to predict the state $\s^* = s(\mupar^*)$
by computing
\begin{equation}
\s^* = \Phi g(\mupar^*).
\end{equation}
This method is called POD with interpolation due to this regression
function acting on the latent variables. Common choices for the
interpolatory map are radial basis functions~\cite{buhmann2003radial}, Gaussian
process~\cite{guo2018reduced, ortali2020mine}, cubic
splines~\cite{lyche1973convergence, gadalla2021les}, or
artificial neural networks~\cite{pichi2021artificial,
  salvador2021non}. We are going to show how to
compute an efficient approximation of such a map by exploiting only the
directions of maximal variations in a multi-fidelity setting, without
the need to perform additional simulations.

\subsection{Parameter space reduction through active subspaces}

Active subspaces~\cite{constantine2015active} is a gradient-based technique for
parameter space reduction~\cite{tezzele2022aroma16}. Let $f: \R^M \to \R$ be a scalar function
of interest. Through the
eigendecomposition of the uncentered
covariance matrix of the gradients, also denoted as the second moment
matrix of $\nabla f$, we can identify the direction of 
maximal variation of $f$ along the parameter space. Let $Q$ be such matrix:
\begin{equation}
Q := \mathbb{E}_{\rho}\, [\nabla_{\mupar} f \, \nabla_{\mupar} f
^T] =\int (\nabla_{\mupar} f) ( \nabla_{\mupar} f )^T\, \rho\  d\mathcal{L},
\end{equation}
where $\mathbb{E}_{\rho}$ stands for the expected values with respect
to the probability density function $\rho$, and $d\mathcal{L}$ is the
Lebesgue measure. The function $\rho$ characterizes the distribution
of the input parameters. $Q$ is symmetric positive definite and can be
decomposed as $Q = \mathbf{W} \Lambda \mathbf{W}^T$. Analogously to
what we have done for the POD, we retain the
first $r^\prime \ll M$ eigenvectors, $\mathbf{W}_{r^\prime}$, and we
use them to project the data onto the so-called active subspace, that
is $\etapar := \mathbf{W}_{r^\prime}^T\mupar \in \R^{r^\prime}$. The
truncation rank $r^\prime$ can be selected a priori or through
the spectral decay of the matrix $Q$. With this projection we are
essentially discarding the directions of the parameter space along
which $f$ is constant or almost constant. Finally we can construct a
ridge approximation $g^\prime$ of the function of interest, that is
\begin{equation}
f (\mupar) \approx g^\prime (\mathbf{W}_{r^\prime}^T\mupar) = g^\prime
(\etapar),
\end{equation}
which lives on a low-dimensional space~\cite{pinkus2015ridge}. Usually
a Gaussian process is used to construct such response surface. In the
next section we are going to show how to exploit $g^\prime$ to build a
low-fidelity model and increase the accuracy of the reduced state
predictions.

AS has been successfully used to approximate the POD modal
coefficients under some assumptions on the size of the
dataset~\cite{demo2019cras}. When in presence of different fidelity
models for $f$ a multi-fidelity model for the computation of the AS
can be used~\cite{lam2020multifidelity}. Other applications in naval engineering
can be found in~\cite{tezzele2018dimension, tezzele2018ecmi}, while
for coupling of AS with model order reduction methods
see~\cite{tezzele2018combined, tezzele2020enhancing}.

\subsection{Nonlinear autoregressive multi-fidelity GP}

The nonlinear autoregressive multi-fidelity
Gaussian process regression (NARGP) scheme was proposed
in~\cite{perdikaris2017nonlinear}. Let us consider the input/output pairs
corresponding to $p$ levels of increasing fidelity, that is
\begin{equation}
  \mathcal{S}_q = \{ x_i^q, y_i^q \}_{i=1}^{N_q} \subset
\mathcal{X} \times\R \subset\R^M \times \R, \qquad \text{ for } q \in \{1, \dots, p \},
\end{equation}
where $y_i^q = f_q (x_i^q)$, and $q=1$ stands for the lowest
fidelity. Let $\pi:\R^M\times\R\rightarrow\R^M$ be the map projecting
the data onto the first $m$ coordinates corresponding to the input parameters. We
 assume the following hierarchical structure:
\begin{equation}
  \pi(\mathcal{S}_p)\subset\pi(\mathcal{S}_{p-1}) \subset \dots
  \subset \pi(\mathcal{S}_1),
\end{equation}
so that we only need a small number of high-fidelity data with respect
to the low-fidelity ones. The hierarchy is due to the autoregressive
nature of the method.
The key step is to assign to each fidelity model $f_{q}$ a Gaussian
process defined by the mean field $m_{q}$, and by the kernel $k_{q}$, as follows:
\RB{\begin{equation}
  y_q(\bar{x})-\epsilon \sim \mathcal{GP}(f_q(\bar{x})|m_q(\bar{x}), k_q(\theta_q))\quad \forall q \in \{1, \dots, p \},
\end{equation}}
where $\epsilon\sim\mathcal{N}(0, \sigma^{2})$ is a noise term and
\begin{equation}
  \bar{x}:=\begin{cases}
    (\mathbf{x},f_{q-1}(\mathbf{x}))\in\R^d\times\R, & q>1\\
    \mathbf{x}\in\R^d, & q=1
  \end{cases}.
\end{equation}

The idea is to build a low-fidelity model exploiting 
parameter space reduction through active
subspaces. We are going to consider only two
fidelity levels, where the lowest one is built without the need of new
simulations coming from simplified models, but instead is a constant
extension along the inactive subspace of the regression built along the
active subspace. We call this approach NARGPAS. We use the high-fidelity data as training
set, following the algorithm presented in~\cite{romor2021multi}, with
the following multi-fidelity model: 
\begin{equation}
g_{\text{MF}}=((f_{H}|x^{H}_{i}, y^{H}_{i}),\
    (f_{L}|x^{L}_{i}))\sim (\mathcal{GP}(f_{H}|m_{H}, \sigma_{H}),
    \mathcal{GP}(f_{L}|m_{L}, \sigma_{L})),
\end{equation}
where the $H$ and $L$ denote the high and low-fidelity,
respectively. The scalar quantity of interest we are going to model are the
reduced state variables. The low-fidelity is built by extending on the whole
parameter space a one-dimensional response surface constructed over
the AS corresponding to each POD coefficient, thus the name POD-NARGPAS. As we are going to show
in the section devoted to the numerical 
results, the new proposed data-driven approach outperforms the more classical
single-fidelity, without the need of any additional
simulation. We remark that POD-GPR, for naval engineering problems, has proven
better than RBF and linear interpolation in~\cite{ortali2020mine}.

For the computation of the active subspace we used ATHENA\footnote{Available at
\url{https://github.com/mathLab/ATHENA}.}~\cite{romor2020athena}. For
the construction of the reduced order models we used
EZyRB\footnote{Available at
  \url{https://github.com/mathLab/EZyRB}.}~\cite{demo2018ezyrb}.

\section{Bayesian optimization}
\label{sec:ironth_bayes}

To minimize the model functions of interest we use Bayesian
optimization~\cite{pelikan1999boa, shahriari2015taking,
  frazier2018bayesian}, which we are going to briefly present in this
section. It is a class of machine-learning-based derivative-free
global optimization methods. One of the main
assumptions is that we do not have any information about the structure
of the function to optimize, so it is intended as a black-box.
Bayesian optimization was first introduced in~\cite{kushner1964new,
  zhilinskas1975single, movckus1975bayesian}, and
successively made popular in~\cite{jones1998efficient} in the context
of efficient global optimization.

Mathematically, we are
considering the problem of finding a global minimizer of an unknown
function $f : \Omega \subset \R^M \to \R$, that is
\begin{equation}
x_{\text{opt}} = \argmin_{x \in \Omega} \, f(x),
\end{equation}
where $\Omega \subset \R^M$ is the design space of interest, which in
our case is the parameter space $\mathcal{P}$. More
general settings can include design spaces with less regularity due to
the presence of possible nonlinear constraints. Here
we consider a sequential search algorithm which selects the next
location where to query $f$. The Bayesian posterior represents the
best current knowledge of the function to optimize. The locations are
selected by evaluating an acquisition function $\alpha : \Omega \to
\R$ which leverages the uncertainty of the posterior to guide the
exploration of the design space.
In \RB{Algorithm~\ref{algo:ironth_byes_opt}} we sketch a 
pseudo-code to highlight the main steps of the whole process. For the
actual implementation we used the one provided by Emukit~\cite{emukit2019}.

\begin{algorithm}[htb]
  \caption{Bayesian optimization pseudo-code.}
  \label{algo:ironth_byes_opt}

  \hspace*{\algorithmicindent} \textbf{Input}:\\
  \hspace*{\algorithmicindent}\hspace*{10pt} model function $f$ to minimize \\
  \hspace*{\algorithmicindent}\hspace*{10pt} acquisition function $\alpha$\\
  \hspace*{\algorithmicindent}\hspace*{10pt} initial number of evaluation points $n_0$ \\
  \hspace*{\algorithmicindent}\hspace*{10pt} maximum number of iterations $N$ \\
  \hspace*{\algorithmicindent} \textbf{Output}:\\ 
  \hspace*{\algorithmicindent}\hspace*{10pt} smallest value for $f(x)$

  \begin{algorithmic}[1]
    \State Compute the Gaussian process prior on $f$.
    \State Evaluate $f$ at the initial $n_0$ points.
    \State Set $n = n_0$.
    \While{$n \leq N$}
    \State Update the posterior probability distribution on $f$ with all the
    available evaluations.
    \State Using the current posterior distribution find $\x_n = \argmax_x \, \alpha(x)$.
    \State Evaluate $y_n = f(x_n)$.
    \State Increment $n$: $n = n+1$.
    \EndWhile
    \State \Return either the point with the smallest $f(x)$, or the one with the
  smallest posterior mean.

  \end{algorithmic}
\end{algorithm}

As acquisition function we adopt the expected improvement (EI), which is
one of the most commonly used. Let
\begin{equation}
f_n^*(x) := \min_{l \leq n} f(x_l),
\end{equation}
be the value of the smallest observed values. We want to perform a new
evaluation, say $y = f(x)$, which has the highest expected improvement
defined as:
\begin{equation}
\alpha (x) := \text{EI}_n(x) := \mathbb{E}_n \left [ (f_n^*(x) - f(x))_+ \right ].
\end{equation}
With $\mathbb{E}_n$ we denote the expectation taken under the
posterior distribution given the first $n$ evaluations, that is
$\mathbb{E}_n [\cdot] = \mathbb{E}_n [\cdot | x_1, \dots, x_n, y_1,
\dots, y_n]$\RB{, while with $(\cdot)_+ = \max (\cdot, 0)$ we denote
  the positive part}. This acquisition function can be computed in
closed form~\cite{jones1998efficient}. The actual next point $x_{n+1}$ to evaluate is
then given by
\begin{equation}
  x_{n+1} = \argmax_x \text{EI}_n(x).
\end{equation}
To find it we exploit the simpler structure of $\alpha$ with respect
to the target function $f$, which allows for inexpensive evaluations
and also easy computation of first and second derivatives.

\begin{figure}[h!]
\centering
\includegraphics[width=.49\textwidth]{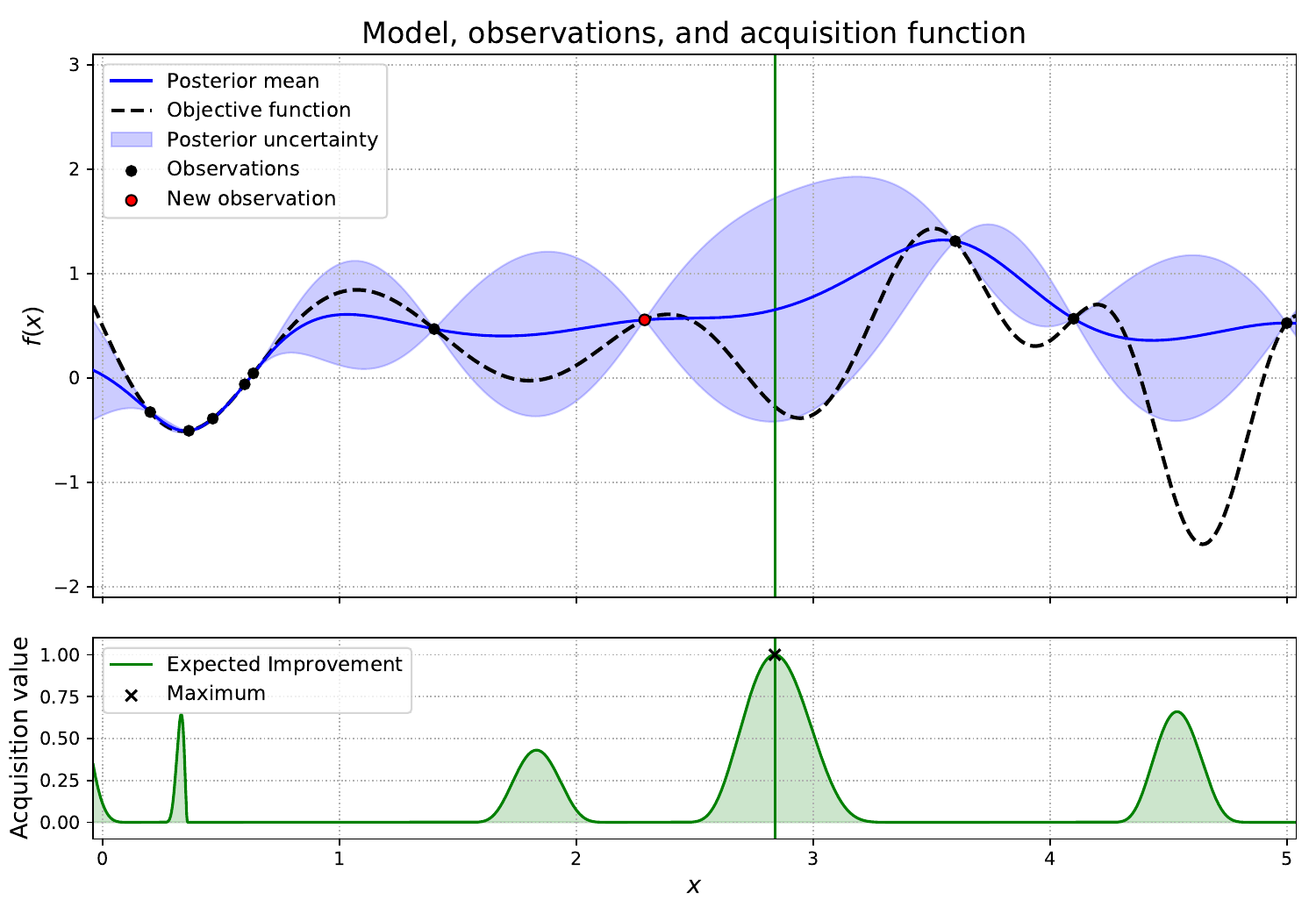}
\includegraphics[width=.49\textwidth]{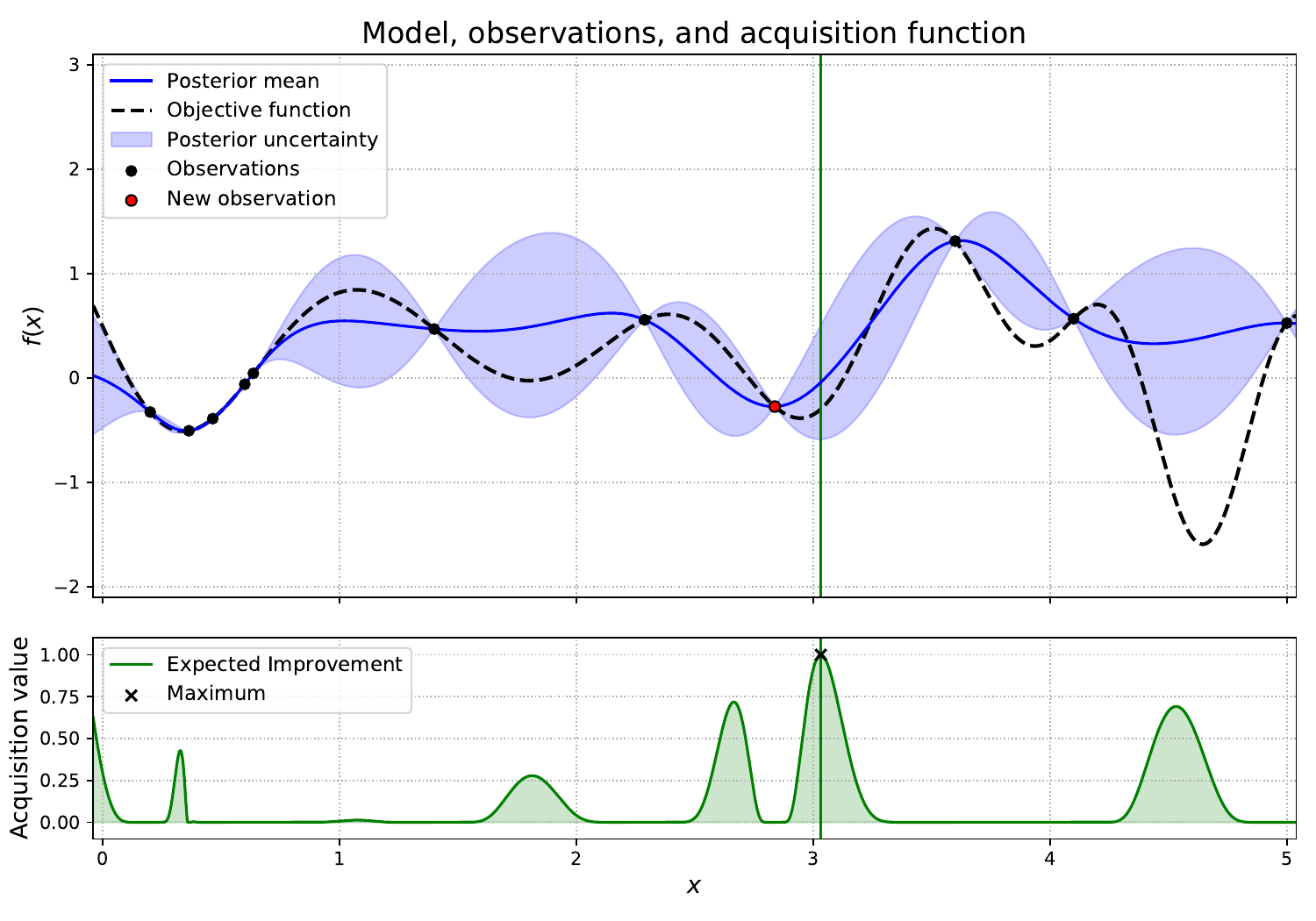}
\caption{\RB{Illustration of the Bayesian optimization procedure for a
    given iteration (on the left) and the next iteration (on the right). The top part of the plots shows the estimated mean and
  confidence intervals of the unknown objective function (in dashed
  line). In the bottom part of the plots we show the acquisition
function (in green) and its maximum.}}
\label{fig:ironth_acq_gpopt}
\end{figure}

In \autoref{fig:ironth_acq_gpopt} we show $2$ iterations of a Bayesian
optimization procedure for a test function, highlighting both the prediction and the
value of the acquisition function. We notice that the acquisition is high where the model
predicts a low objective (so-called exploitation) and where the prediction uncertainty
is high (so-called exploration). Note that the area on the far left remains unsampled,
as while it has high uncertainty, it is predicted to offer a smaller
improvement over the best observation. 

\section{Numerical results}
\label{sec:ironth_results}

In this section we are going to apply the optimization pipeline
described in \autoref{sec:ironth_opt_pipeline} to a midship section and
to the parametrized hull depicted in \autoref{fig:ironth_c6298}. The
first test case is intended to be illustrative due to its smaller size
and the possibility of showing all the error comparisons. We
considered a midship section of a simplified hull with $20$
parameters. It presents all the challenges of an entire hull in terms
of stress distribution but with less degrees of freedom, allowing for
faster high-fidelity simulation. The second test case is a real
passenger ship hull parametrized with $16$ parameters, keeping the
number of parameters comparable to the previous test. We perform a
discrete mono-objective optimization of the mass of the parametrized
regions, considering stability constraints.

\subsection{Objective function definition}
The target function we are going to minimize is the total mass $m$ of the
parametrized regions of the hull plus the mass of the buckling
stiffeners needed to stabilize the buckled elements of the entire
hull. Other possible choices are available within the optimization
framework such as the deflection at a selected point for a given loading condition, and
the total cost related to the parametrized decks. The latter considers both the
acquisition cost of the metal raw materials and the manufacturing cost
for the installation of the steel plates and buckling stiffeners, specific for each
shipyard. For industrial reasons in this work we are not going to
present the results for the cost optimization, but we want to emphasize that the
framework is very versatile and allows the use of different target functions.

As for the stability constraints we set some thresholds for the Cauchy stress
tensor components in order to count how many plate elements are
yielded for a prescribed set of loading conditions, which in the
present work are hogging and sagging. See \autoref{fig:ironth_c6298_vm} for an example of
hogging condition. Given the symmetric Cauchy stress tensor in the global
reference frame, whose components for a 
single element are 
\begin{equation}
\sigma =
\begin{bmatrix}
\sigma_x & \tau_{xy} & \tau_{xz} \\
\tau_{xy} & \sigma_y & \tau_{yz} \\
\tau_{xz} & \tau_{yz} & \sigma_z ,
\end{bmatrix},
\end{equation}
we define an element yielded if at least one of the following
conditions is not satisfied:
\begin{align}
  -245 &\leq \sigma_i \leq 245, \qquad \text{ for } i \in \{x, y, z\},\\
  -153 &\leq \tau_i \leq 153, \qquad \text{ for } i \in \{xy, xz, yz\}, \\
\sigma_{\text{VM}} &:= \sqrt{\sigma_x^2 + \sigma_y^2 - \sigma_x
         \sigma_y + 3 \tau_{xy}^2} \leq 307,
\end{align}
where $\sigma_{\text{VM}}$ stands for the von Mises yield
criterion. The thresholds above are characteristics of the high strength structural steel.
The actual constraint for the optimization is the maximum number
$N_{\text{max}}^y$ of elements that can yield. Exactly the same
is done for the buckling usage
factors associated to each element. An element is considered buckled,
based on the DNV GL classification rules\footnote{Det Norske Veritas (DNV) Rules for
  Ships, part 3, chapter 1, section 13: Buckling control.}, if at least one of the $11$
components of the buckling usage factors tensor is greater than $1$, for at
least one loading condition. Such tensor is computed as a function of
the Cauchy stress tensor. The maximum number of allowed buckled
elements is $N_{\text{max}}^b$.  To incorporate these stability
constraints we penalize the objective function $f_{\text{obj}}$ with a parabolic function
depending \RB{on the violated} constraint. Its expression is the following
\begin{equation}
  \label{eq:ironth_obj_f}
f_{\text{obj}} (\mupar) := m(\mupar) + m_{\text{bs}} N^b(\mupar) + c_y (N^y(\mupar) - N_{\text{max}}^y)_+^2
+ c_b (N^b(\mupar) - N_{\text{max}}^b)_+^2, 
\end{equation}
where $m(\mupar)$ is the mass of the parametric decks, $N^y (\cdot)$
and $N^b (\cdot)$ denote the number of yielded and buckled 
elements, respectively, $m_{\text{bs}} = 2.968$~\si{kg} is the mass of a single
buckling stiffener, $c_y = 1$, $c_b = 0.001$, and $(\cdot)_+ =
\max (\cdot, 0)$ stands for the positive part. The coefficients $c_y$
and $c_b$ are prescribed by the user depending on the order of
magnitude of the other terms. With this formulation we allow local stress peaks
to exceed the rule-based strength limits because we are able to
account for the additional mass needed to stabilize the affected
elements.

\subsection{Midship section}
\begin{figure}[htb]
    \centering 
    \includegraphics[trim=0 0 0 0, clip, width=.49\textwidth]{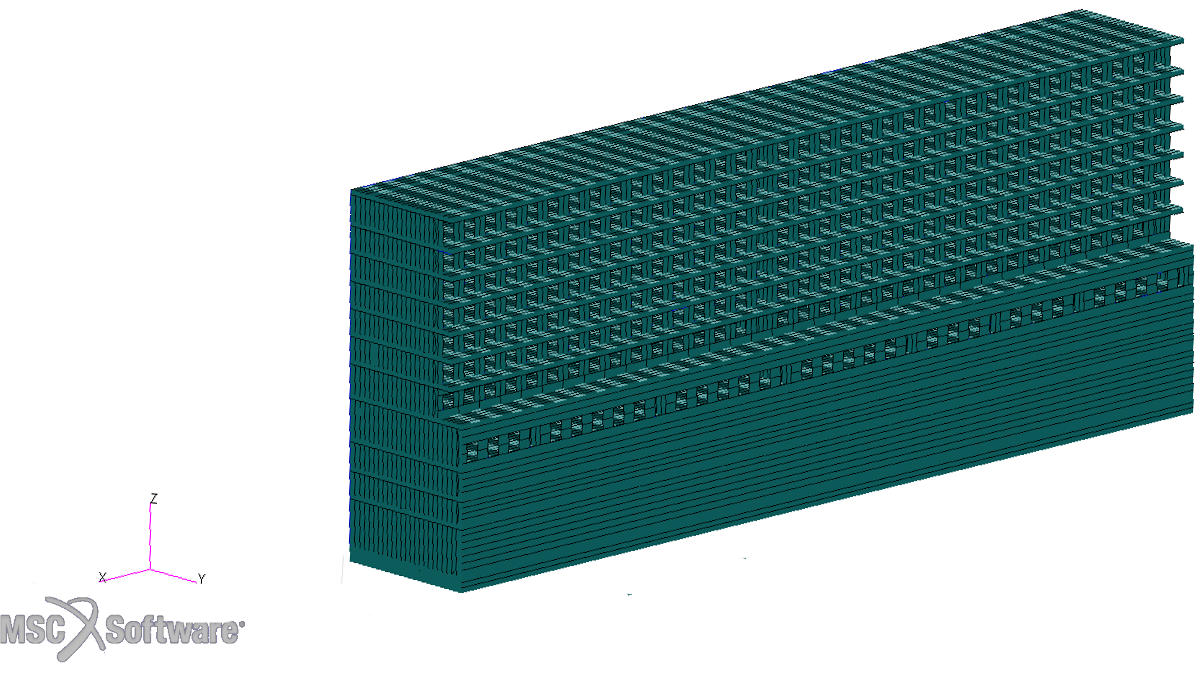}\hfill
    \includegraphics[trim=0 0 0 0, clip, width=.49\textwidth]{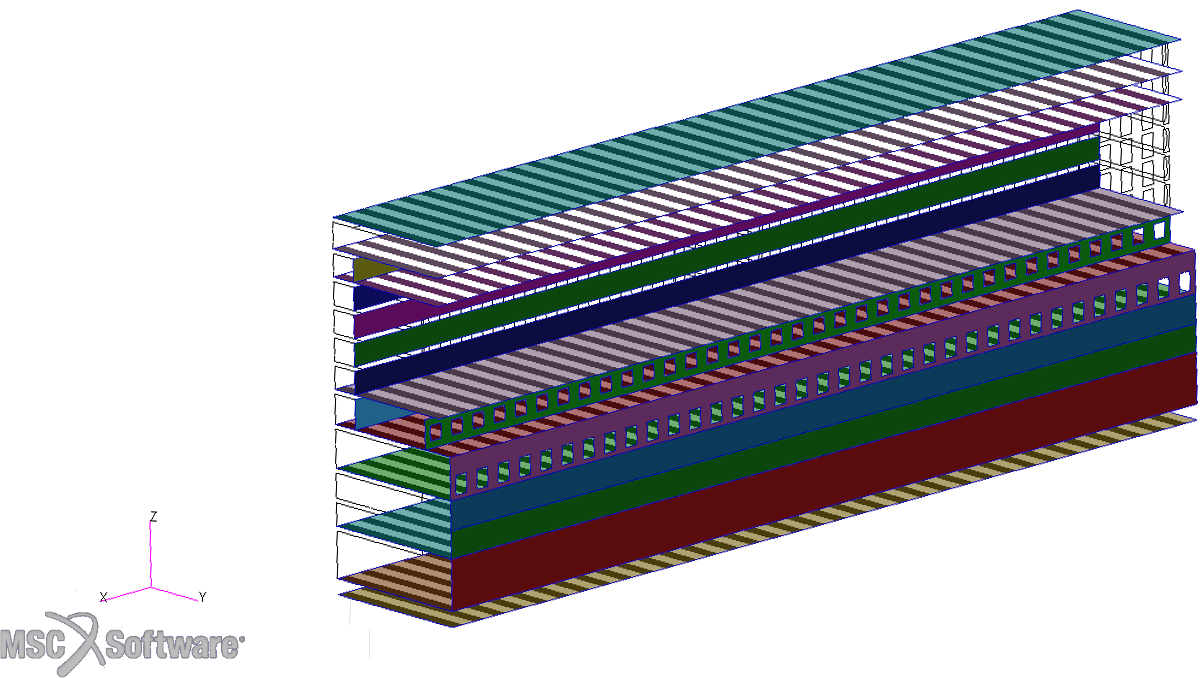}
    \caption{A complete view of the midship section on the left, and the highlight of the 20 parametrized elements groups on the right.}
    \label{fig:ironth_midship}
\end{figure}

For the testing and tuning phases of the pipeline development, we
consider a smaller model to obtain faster high-fidelity
evaluations. This midship section, shown in the left panel of
\autoref{fig:ironth_midship}, approximates a quarter of the
ship's elements, excluding the bow and stern segments; the finite
element problem applies mirroring along both negative $x$ and negative
$y$ coordinates, as the origin represent the ship's center. Due to
this simplification and the very regular internal structure, this
midship section contains about $\frac{1}{10}$ of the elements of the
full ship and the required high-fidelity solve time is about
$\frac{1}{20}$.

\begin{figure}[htb]
    \centering
    \includegraphics[trim=0 0 0 0, clip, width=.49\textwidth]{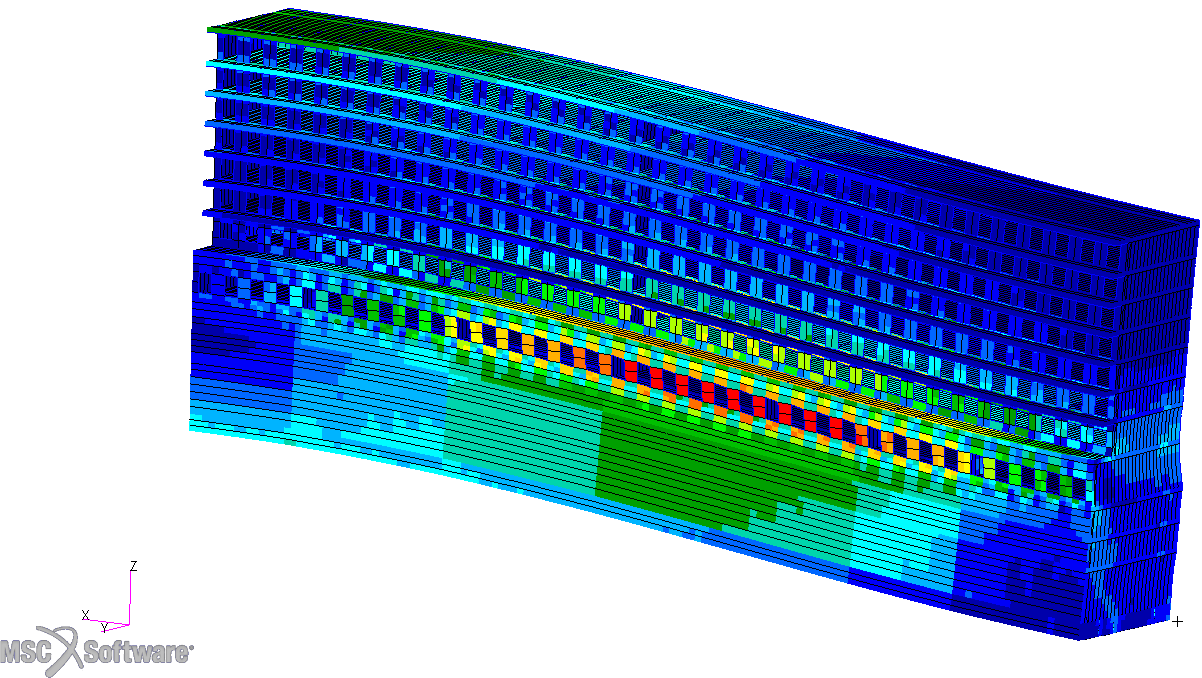}\hfill
    \includegraphics[trim=0 0 0 0, clip, width=.49\textwidth]{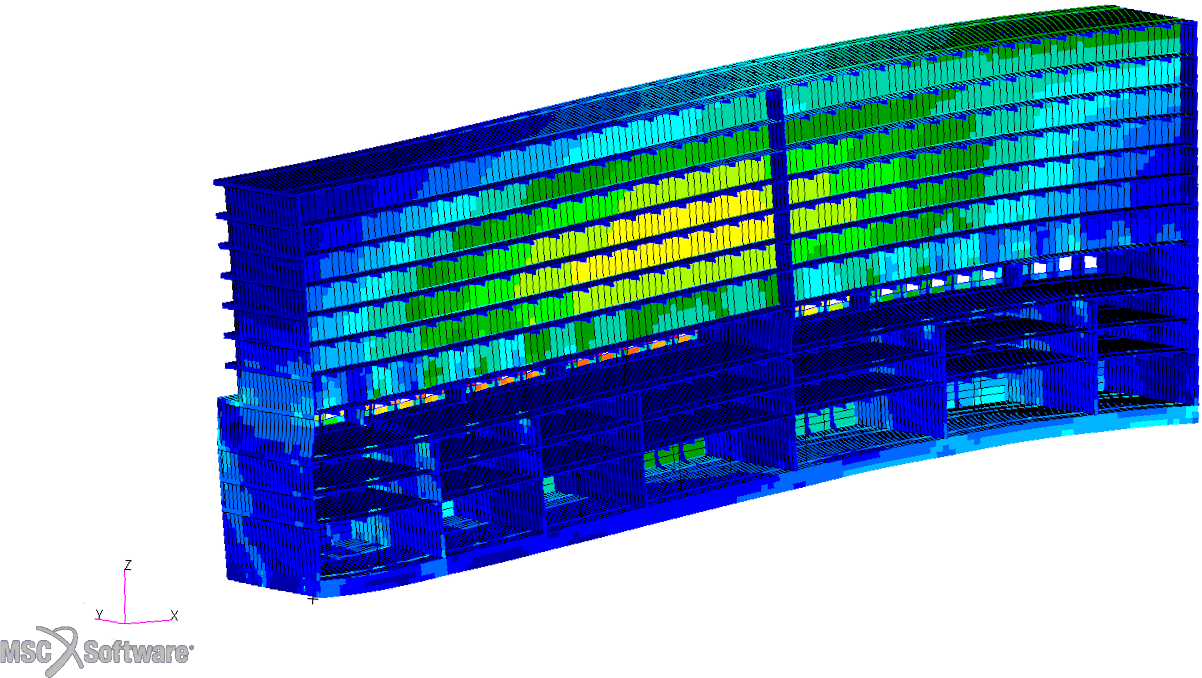}
    \caption{The external view of the sagging midship section on
          the left, and the internal view of the same configuration on
          the right. Displacements are magnified. Colors refer to the
          von Mises criterion.} 
    \label{fig:ironth_midship_sag}
\end{figure}

The parametrized sections, depicted in the right panel of \autoref{fig:ironth_midship}, were chosen with the aim of studying the
structural resilience against longitudinal loads: all groups of
elements have the largest span along the $x$ direction. Of these, the
$9$ decks included are equally distributed between the lowermost,
uppermost and central; whereas the groups which span vertically
comprise the $4$ external plankings between decks $1$ and $5$, the $6$
internal twin decks between decks $5$ and $11$, and finally the external
twin deck between decks $5$ and $6$. \autoref{tab:ironth_20_params_hull}
describes the actual parameter space $\mathcal{P} \subset \R^{20}$ used, the location
of the area affected by a specific parameter and its default value. 

\begin{table}[htb]
    \centering
    \caption{Parameters description of the hull for the midship test case. All data are in \si{mm}.\label{tab:ironth_20_params_hull}}
    \begin{tabular}{ c c c c c }
        \hline
        \hline
        Parameter & Region & Default thickness & Lower bound & Upper bound \\
        \hline
        \hline
        \rowcolor{Gray}
        $\mupar_{1}$ & Deck $12$ & 6.0 & 5.0 & 20.0 \\
        $\mupar_{2}$ & Deck $11$ & 6.0 & 5.0 & 20.0 \\
        \rowcolor{Gray}
        $\mupar_{3}$ & Deck $10$ & 5.5 & 5.0 & 15.0 \\
        $\mupar_{4}$ & Twin deck $5$,$6$ & 5.0 & 5.0 & 15.0 \\
        \rowcolor{Gray}
        $\mupar_{5}$ & Twin deck $6$,$7$ & 5.0 & 5.0 & 15.0 \\
        $\mupar_{6}$ & Twin deck $7$,$8$ & 5.0 & 5.0 & 15.0 \\
        \rowcolor{Gray}
        $\mupar_{7}$ & Bottom & 14.0 & 12.0 & 20.0 \\
        $\mupar_{8}$ & Deck $1$ & 13.0 & 12.0 & 20.0 \\
        \rowcolor{Gray}
        $\mupar_{9}$ & Deck $2$ & 6.0 & 5.0 & 15.0 \\
        $\mupar_{10}$ & Twin deck $8$,$9$ & 5.0 & 5.0 & 15.0 \\
        \rowcolor{Gray}
        $\mupar_{11}$ & Twin deck $9$,$10$ & 5.0 & 5.0 & 15.0 \\
        $\mupar_{12}$ & Twin deck $10$,$11$ & 5.0 & 5.0 & 15.0 \\
        \rowcolor{Gray}
        $\mupar_{13}$ & Deck $4$ & 5.0 & 5.0 & 15.0 \\
        $\mupar_{14}$ & Deck $5$ & 5.0 & 5.0 & 15.0 \\
        \rowcolor{Gray}
        $\mupar_{15}$ & Deck $6$ & 5.0 & 5.0 & 15.0 \\
        $\mupar_{16}$ & External twin deck $5$,$6$ & 8.0 & 5.0 & 15.0 \\
        \rowcolor{Gray}
        $\mupar_{17}$ & External planking $4$,$5$ & 10.0 & 8.0 & 15.0 \\
        $\mupar_{18}$ & External planking $3$,$4$ & 10.0 & 8.0 & 15.0 \\
        \rowcolor{Gray}
        $\mupar_{19}$ & External planking $2$,$3$ & 12.0 & 8.0 & 15.0 \\
        $\mupar_{20}$ & External planking $1$,$2$ & 12.0 & 8.0 & 15.0 \\
        \hline
        \hline
    \end{tabular}
\end{table}

In \autoref{fig:ironth_midship_sag} the von Mises stress values are shown
under sagging. The effect of shear stresses on the external plankings
and twin decks is particularly pronounced on the central twin decks,
as well as the normal stresses on the central and uppermost
decks. These images were obtained for the default parameters
configuration, where the mentioned groups are very thin.

\begin{figure}[htb]
\centering
\includegraphics[width=.7\textwidth]{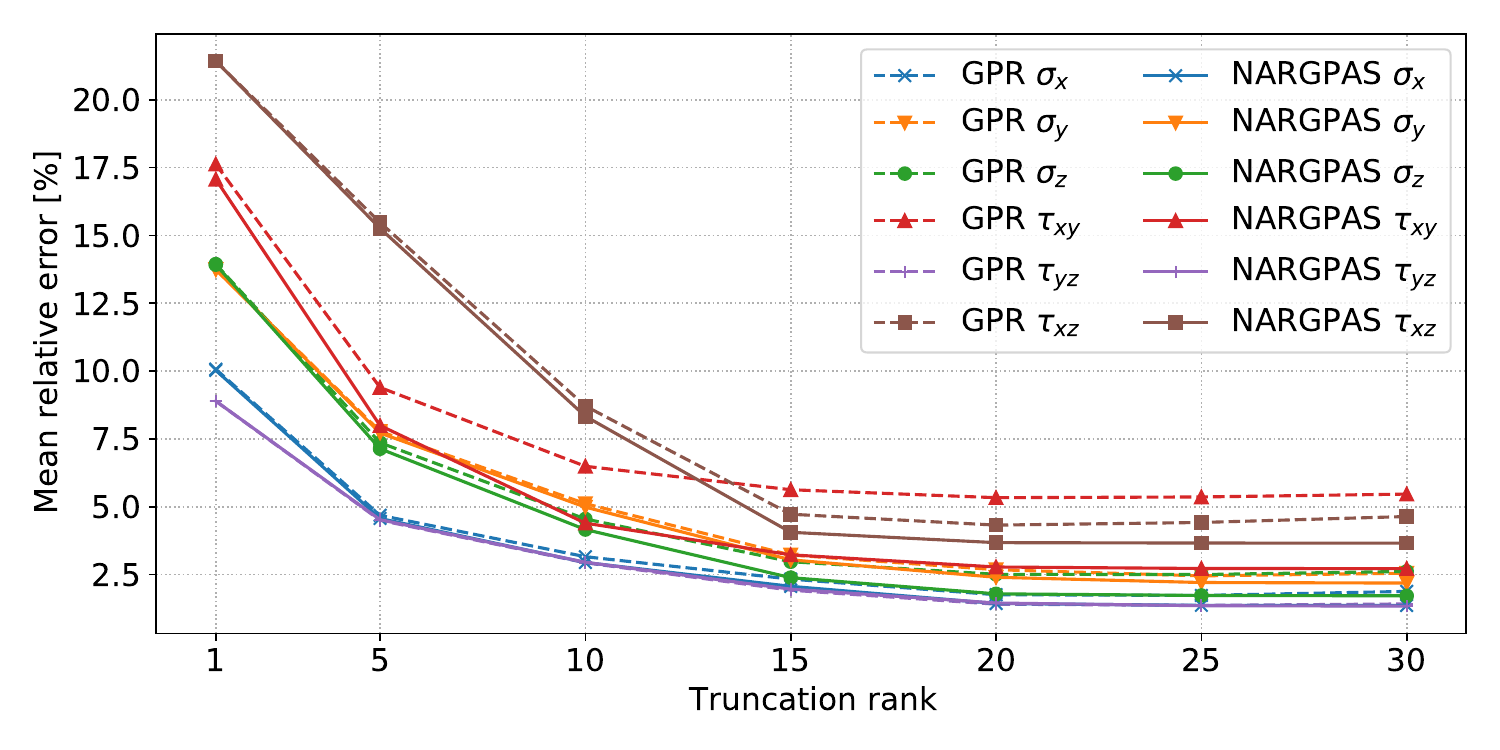}
\caption{The mean relative $L^2$ error of the stress field varying the truncation rank, for 
GPR (dashed lines) and NARGPAS (continuous lines) interpolators. For each feature, the values are obtained by averaging 
the test scores from a $5$-fold cross-validation experiment.}
\label{fig:ironth_midship_truncation_ranks}
\end{figure}

We consider an initial database of $300$ high-fidelity solutions, and
we adopt a $5$-fold cross validation scheme in order to estimate the
prediction errors for  the stress tensor components. In
\autoref{fig:ironth_midship_truncation_ranks} we compare  the mean
relative $L^2$ errors of each component, obtained with GPR and NARGPAS 
interpolators for different POD truncation ranks. The general trend is a reduction 
of the errors as the truncation rank increases, with NARGPAS achieving better error 
values than GPR especially for $\tau_{xy}$ with a $50\%$ reduction.

\begin{figure}[htb]
\centering
\includegraphics[width=.49\textwidth]{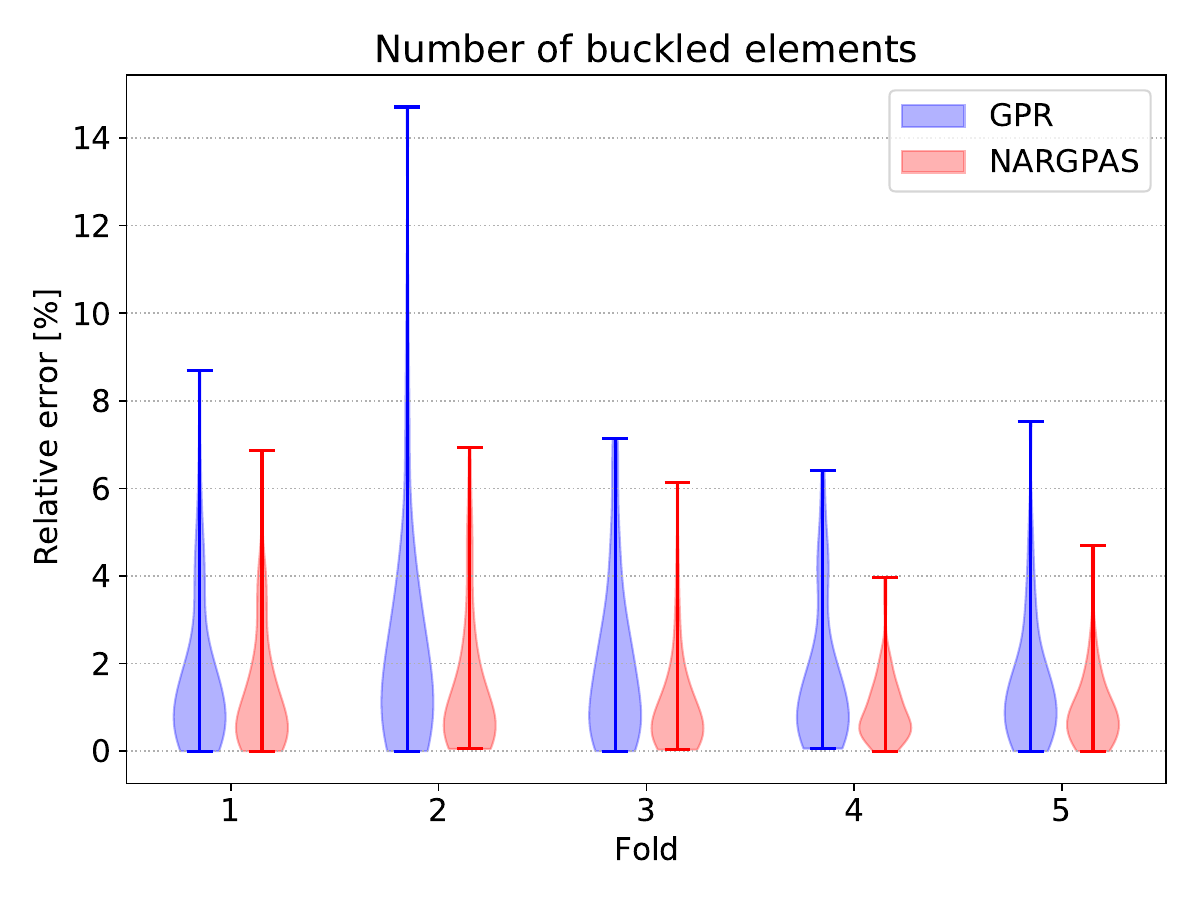}
\caption{Violin plot describing the distribution of the relative error
  in the prediction of the number of buckled elements. Data are
  divided by fold and methods used.}
\label{fig:ironth_midship_buckled_error}
\end{figure}

In \autoref{fig:ironth_midship_buckled_error}, 
we compare the distribution of the relative prediction errors on the number of 
buckled elements in each fold of the cross validation experiment. As expected from 
the stress tensor components prediction errors, NARGPAS outperforms GPR in every fold 
and achieves both lower maxima and lower spread, presenting a greater
accumulation of the errors near $0$. The POD-NARGPAS approach in general provides
better perfomance both in the approximation of the stress tensor
field and to the derived quantities of interest related to the
stability of the hull.
The results regarding the number of yielded elements are not reported since already with the POD-GPR
method the predictions are very close to the actual values and the
differences are negligible.

After constructing a reduced order model for every component of the
stress tensor we can perform the Bayesian optimization described above. 
\RB{The truncation rank for the stress components is $21$, chosen by analyzing 
  the singular values decay of the snapshot matrix and selecting
  the rank $r$ which manifests the last sharpest decrease between $r$ and $r+1$, among all 
  the stress tensor components. The chosen rank results in a
  cumulative energy greater than $99.99$\%.}  
We emphasize that
every parametric term in \autoref{eq:ironth_obj_f}, except $m(\mupar)$,
is predicted using the precomputed reduced order models. In this way a
single parameter evaluation takes approximately $1$ second. The term
$m(\mupar)$ is computed exactly. For this case we set
$N_{\text{max}}^y = 20$ and $N_{\text{max}}^b = 4000$.

\begin{figure}[h!]
\centering
\includegraphics[width=1.\textwidth]{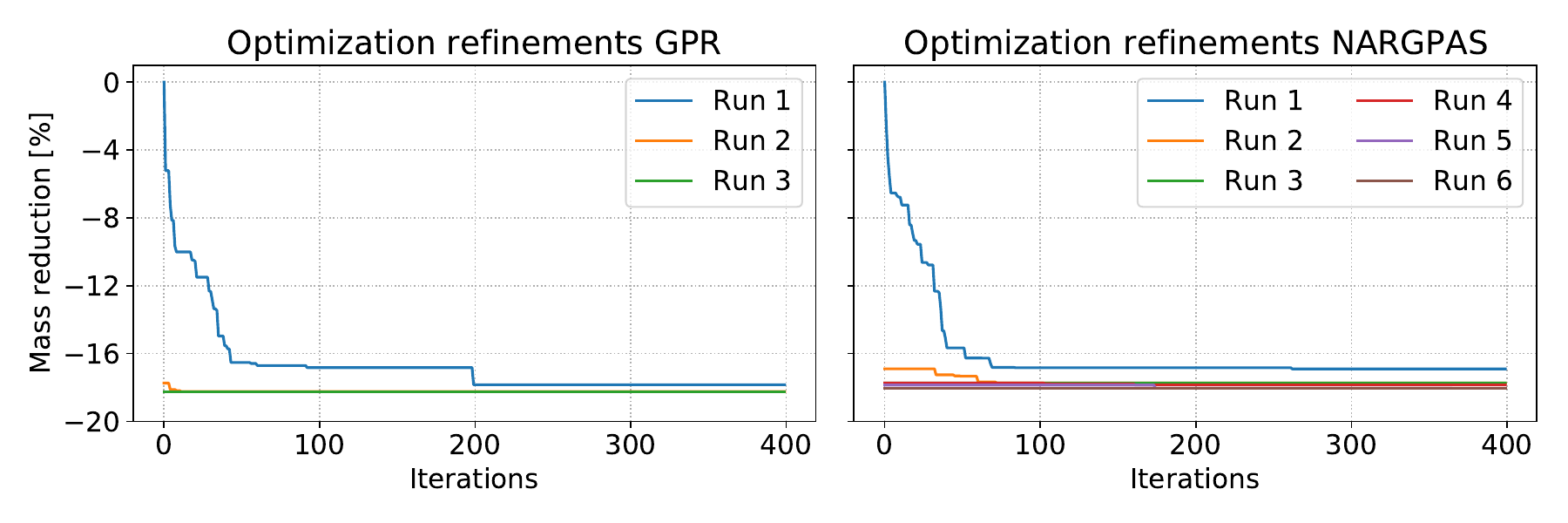}
\caption{Different optimizations runs for the midship test case. On
  the left panel we used the POD-GPR approach, while on the right
  panel the POD-NARGPAS one. The relative reduction is with
  respect to the best sample among the initial solutions database.}
\label{fig:ironth_opt_mass_midship}
\end{figure}

We set the computational budget for the Bayesian optimization to $400$
iterations. We remark that at every iteration the GPR for the target function has to be
recomputed. This means that the GPR construction time
increases at every iteration. In \autoref{fig:ironth_opt_mass_midship} the results for all
the successive optimization runs are depicted. 
After an optimization cycle is completed, a subset of the best low-fidelity configurations evaluated is sent to 
the high-fidelity solver and the results are added to the snapshots database. 
The reduced order models are then updated. Thus the accuracy in the 
neighborhood of the current optimum is increased and the successive
runs will exploit such information by focusing in that specific
region, or will explore different areas if the error committed at the
current optimum is too high. For this test case at every enrichment we
add $20$ high-fidelity evaluations. 
We continue to perform optimization runs followed by the enrichment
phase until the optimizer is not able to find a better point with
respect to the previous run. The main source of errors is in
counting the buckling stiffeners needed to stabilize the buckled
elements as can be observed looking at \autoref{eq:ironth_obj_f} 
and at the previous plots. 

The hull configurations found by the algorithm using GPR and NARGPAS,
although with small differences, modified the default hull in the same
way. In both, the external planking was consistently slimmed down,
together with all decks except the central ones and the uppermost,
which were instead left unchanged. All twin decks were slightly bulked
up, except for the external one which was slimmed down. The final mass
for the configuration found by the GPR interpolator is only $0.23\%$
lower than that obtained with NARGPAS, as can be seen in
\autoref{fig:ironth_opt_mass_midship}. Since both configurations
acted on the default parameter values in the same way, the algorithm
appears robust.

\begin{figure}[h!]
\centering
\includegraphics[width=.49\textwidth]{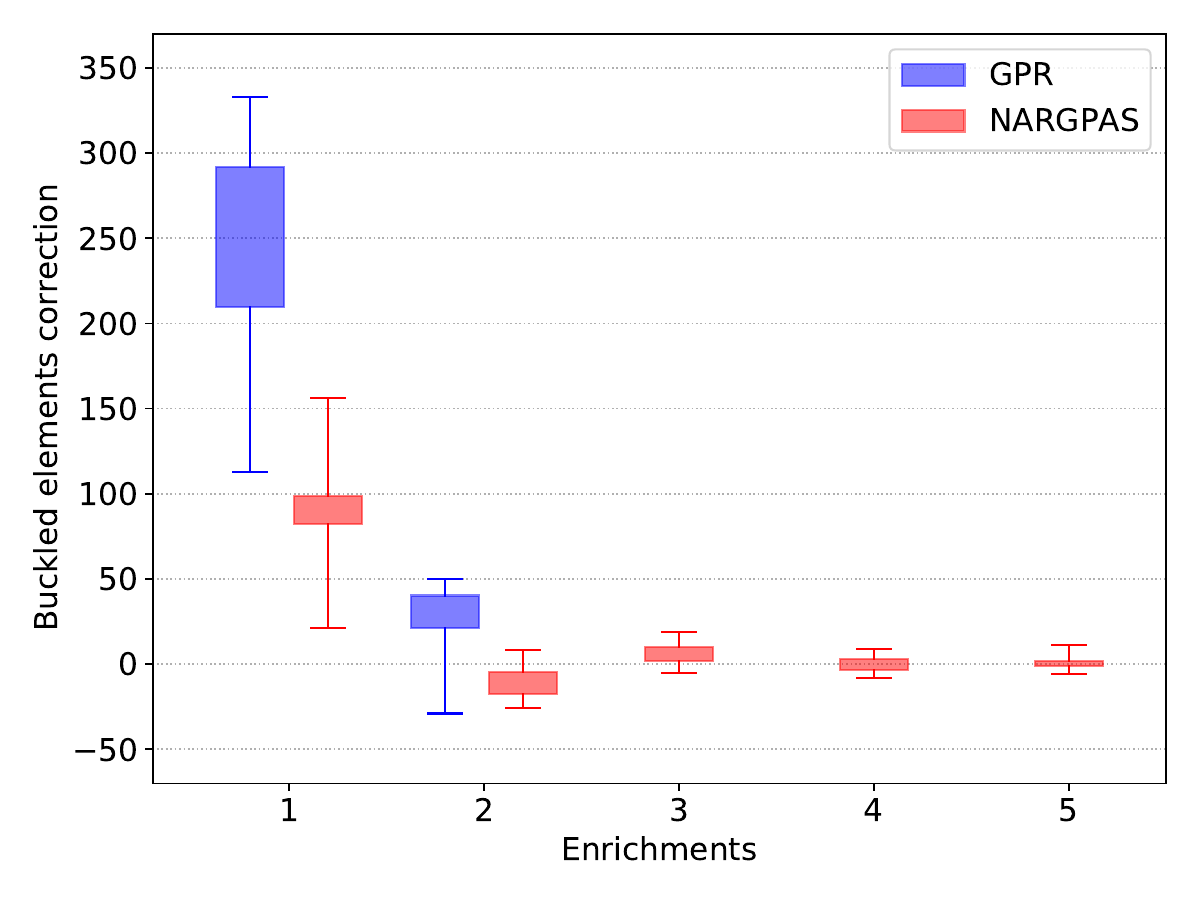}
\caption{Evolution of the prediction error distribution of the number of buckled 
    elements for the validated best candidates, after every optimization-enrichment loop.}
\label{fig:ironth_opt_mass_errorboxes_midship}
\end{figure}

The optimization of GPR-based surrogates ended with the third iteration, while the 
NARGPAS-based ones required six runs. In \autoref{fig:ironth_opt_mass_errorboxes_midship} 
we summarize the evolution of the prediction error distribution on the optimal hulls during 
the optimization-enrichment loop. After each round of Bayesian optimization, 
the best $20$ hull configurations are selected for the enrichment. For each one, 
the difference between the number of buckled elements predicted by the
ROMs used in the previous optimization run and its 
high-fidelity counterpart are stored and the resulting distributions form the 
boxes and whiskers. For both interpolators, the error is largest at the first 
iteration and decreases as the high fidelity database is enriched. GPR-based 
surrogates exhibit larger errors compared to NARGPAS, which despite having required 
a larger number of iterations, seem to provide a faster convergence of the error.

\subsection{Complete hull}

\begin{figure}[htb]
\centering
\includegraphics[trim=0 0 90 40, clip, width=.49\textwidth]{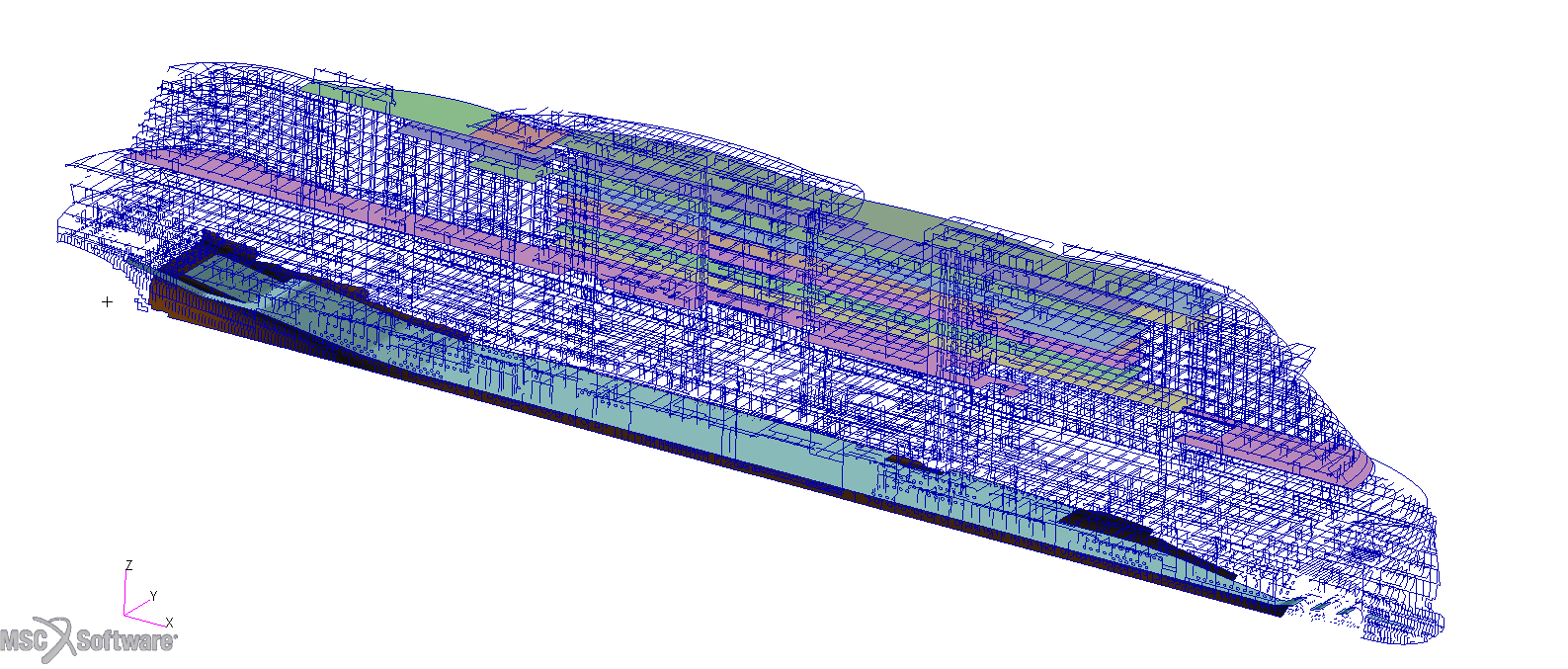}
\caption{View of the hull parametric decks. In solid colors the regions
  of the hull affected by the parameters described in
  \autoref{tab:ironth_16_params_hull}.}
\label{fig:ironth_c6298_params}
\end{figure}

For the second test case we decreased the number of input design
parameters to $16$, but used a larger and more complex model. 
The parametric regions are depicted with solid
colors in the right panel of \autoref{fig:ironth_c6298_params}, while a
deeper description of the range of variations of each parameter
$\mupar_i$, $i=1, \dots, 16$, can be found in
\autoref{tab:ironth_16_params_hull}.
We use the sampling strategy described in \autoref{sec:podi} to generate $300$
samples. \RB{The truncation rank for the stress components is set to $17$,
  for which the cumulative energy of the singular values is greater than $99.99$\%.}. We also
increase the computational budget for the Bayesian optimization to
$600$ evaluations. During the enrichment phase we evaluate the best
$4$ hull configurations found by the optimizer.

\begin{table}[htb]
\centering
\caption{Parameters description for the entire hull. Repeated regions
  correspond to different parts. All data are in \si{mm}.\label{tab:ironth_16_params_hull}}
\begin{tabular}{ c c c c c }
\hline
\hline
Parameter & Region & Default thickness & Lower bound & Upper bound \\
\hline
\hline
\rowcolor{Gray}
$\mupar_1$  & Deck $15$  & 7.5 & 5.0 & 15.0 \\
$\mupar_2$  & Deck $16$  & 8.0 & 5.0 & 20.0 \\
\rowcolor{Gray}
$\mupar_3$  & Deck $17$  & 9.0 & 5.0 & 20.0 \\
$\mupar_4$  & Deck $14$  & 7.5 & 5.0 & 15.0 \\
\rowcolor{Gray}
$\mupar_5$  & Deck $13$  & 7.0 & 5.0 & 15.0 \\
$\mupar_6$  & Deck $12$  & 6.5 & 5.0 & 15.0 \\
\rowcolor{Gray}
$\mupar_7$  & Deck $11$  & 6.0 & 5.0 & 15.0 \\
$\mupar_8$  & Deck $10$  & 5.5 & 5.0 & 15.0 \\
\rowcolor{Gray}
$\mupar_9$  & Deck $17$  & 6.0 & 5.0 & 20.0 \\
$\mupar_{10}$  & Deck $17$  & 15.0 & 5.0 & 20.0 \\
\rowcolor{Gray}
$\mupar_{11}$  & Deck $17$  & 6.0 & 5.0 & 20.0 \\
$\mupar_{12}$  & Deck $16$  & 6.0 & 5.0 & 20.0 \\
\rowcolor{Gray}
$\mupar_{13}$  & Deck $16$  & 6.0 & 5.0 & 20.0 \\
$\mupar_{14}$  & Deck $09$  & 8.0 & 5.0 & 15.0 \\
\rowcolor{Gray}
$\mupar_{15}$  & Deck $01$  & 16.0 & 12.0 & 25.0 \\
$\mupar_{16}$  & Deck $00$  & 20.0 & 12.0 & 25.0 \\
\hline
\hline
\end{tabular}
\end{table}

\begin{figure}[htb]
\centering
\includegraphics[width=1.\textwidth]{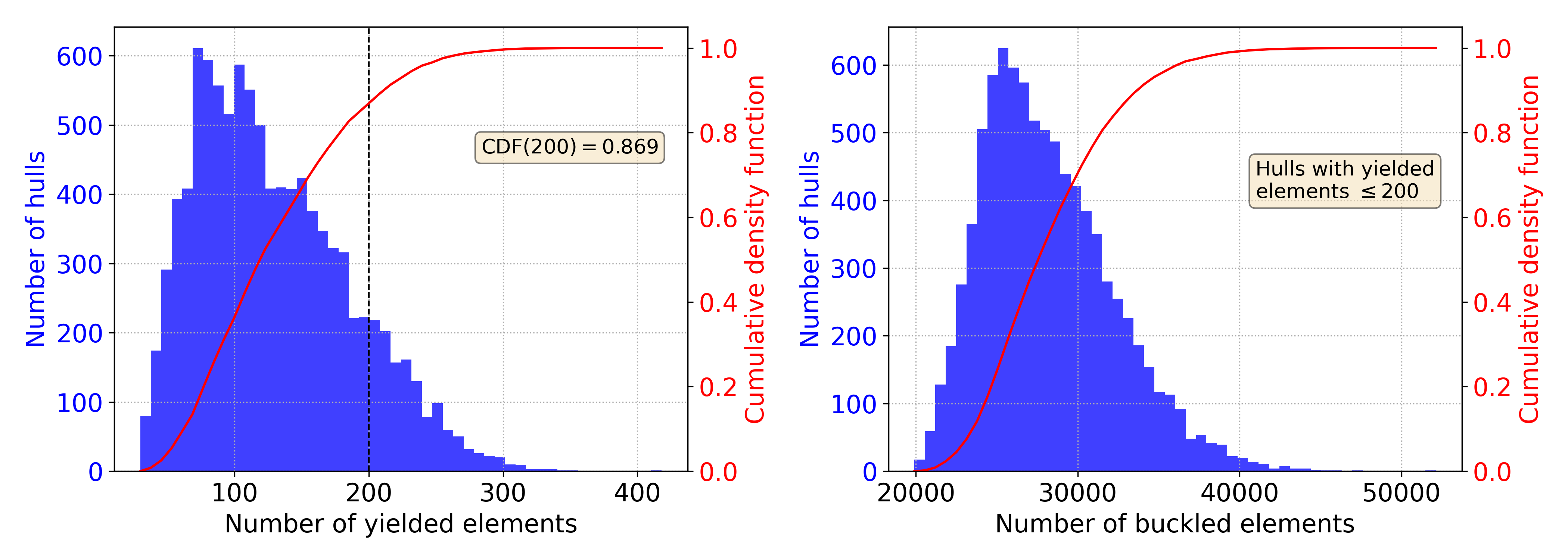}
\caption{Distribution of the hulls with respect to the stability
  constraints. In the left panel we consider the number of yielded elements, while
in the right panel we have the distribution of the parametric hulls
satisfying the first constraint with respect to the number of buckled elements.}
\label{fig:ironth_histograms_constraints}
\end{figure}

To have a better idea on how to set the stability constraints for a
given type of hull, we can plot the distribution of the parametric
hulls with respect to the number of yielded elements (left panel of
\autoref{fig:ironth_histograms_constraints}). We see that if we set
$N_{\text{max}}^y = 200$, the valid samples will be roughly $87\%$ of
all the manufacturable hulls. In the right panel of
\autoref{fig:ironth_histograms_constraints} we plot the distribution of
the hulls satisfying such constraint with respect to the number of buckled elements and we set
$N_{\text{max}}^b = 35000$ accordingly.

\begin{figure}[h!]
\centering
\includegraphics[width=.85\textwidth]{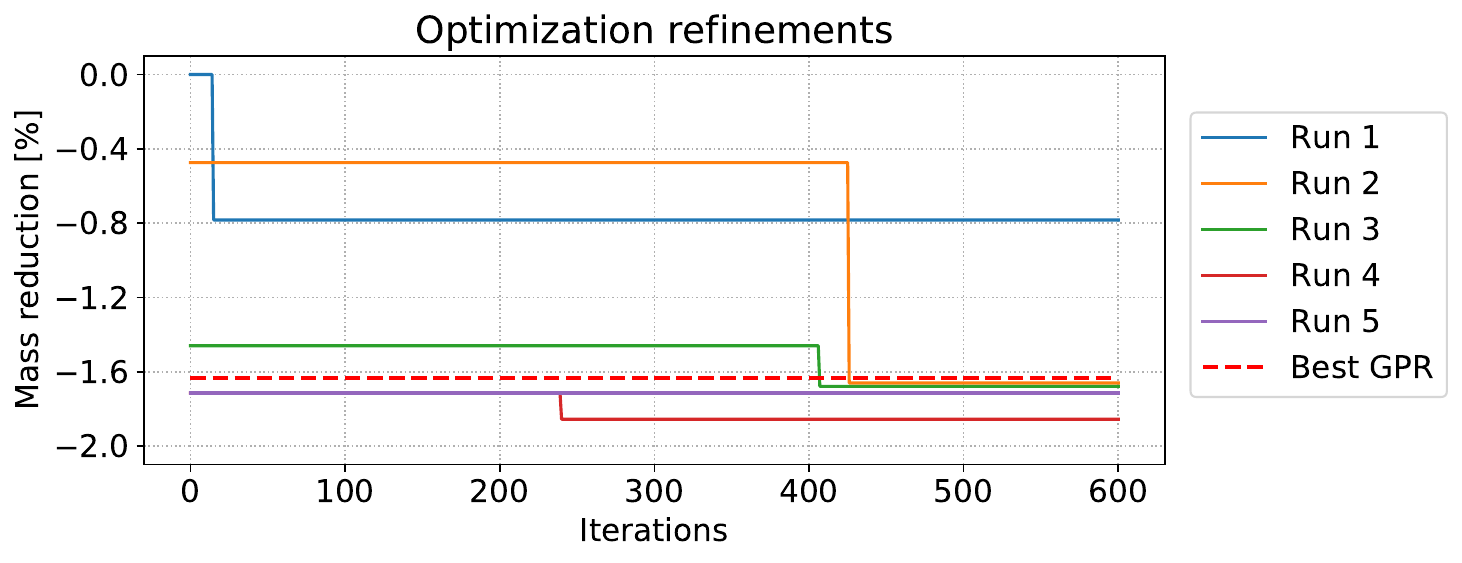}
\caption{Different optimizations runs for the parametrized hull test case. The relative reduction is with
  respect to the best sample among the initial solutions
  database.}
\label{fig:ironth_opt_mass_16_params}
\end{figure}

\autoref{fig:ironth_opt_mass_16_params} depicts the successive runs
performed using the NARGPAS method, compared with the best value found
using GPR to predict the reduced state variables. We need $5$
optimization-enrichment loops to achieve global convergence in the approximated
solution manifold. As we can see the fifth run does not produce any
new improvement. In particular from the figure we conclude that the
third run produces the actual optimum. The slight difference is due to
the approximation error in the computation of the number of buckled
elements, which is not present in the successive run since that
snapshot has been added to the database.

The configurations found by the two optimization loops differ more than in the midship section
case. Both agree on choosing much thinner bottom decks and bulking up
the mid deck. However, the upper decks, which are divided into multiple
sections, show disagreement between the two with how the default
parameter values are changed. Nonetheless, the final mass scores only
differ by $0.08\%$. Absolute values are not reported for industrial
reasons. Even if the mass reduction seems not as remarkable, 
this is due to the fact that the reference hull, present in the initial database, was
obtained after several iterations of the design office and actually
manufactured. With this test case we can see how our framework is able
to find better designs while preserving the stability constraints in a
non-intrusive fashion. 
From the designers' point of view, the changes on the lower and mid
decks are in line with the professional experience accumulated in
years of work spent hand-tuning the initial default parameter
choices. The application of said experience, however, still requires
the designer to proceed with trial and error, which translates to many
days of work due to the long computation time of the high-fidelity
solver and the results analysis. The availability of an automated
procedure which can bootstrap this process in mere hours, and requires
minimal input from the user, speeds up the design phase by a large
margin.

\section{Conclusions and perspectives}
\label{sec:ironth_conclusions}
In this work we proposed a modular data-driven non-intrusive structural optimization
framework for modern passenger ships. We exploit several reduced order
models coupled with parameter space reduction in a multi-fidelity
setting. This new approach is called POD-NARGPAS. We
demonstrated its performance against the more classical POD-GPR.
Our efficient numerical pipeline allows for different discrete
mono-objective optimization given a precomputed set of high-fidelity
simulations. We parametrized the thickness of various regions of the
reference hull and we perform a mass minimization considering
stability constraints such as the total number of yielded and buckled
elements. We allow the presence of unstable elements since we are able
to account for all the necessary interventions to stabilize
them during the optimization process. We provided a comprehensive error analisys for a midship
section test case which presents all the challenges of a complete hull while
keeping low the time to solution. Finally we tested the entire
framework on a real cruise ship and show how the tool is able to find
previously unconsidered designs.

Future works will focus on improving the accuracy of constraints
evaluations, for example with a multi-fidelity approximation of the
scalar output and not only for the reconstruction of the entire
field~\cite{romor2021multi}. Another possibility is the exploitation of local information
with local active subspaces~\cite{romor2021las} or \RA{nonlinear
techniques, based on kernels~\cite{romor2020kas} or
level-sets~\cite{bridges2019active, zhang2019learning},} to further
improve the regression performance of the low-fidelity model. Other physical constraints can
also be considered such as the position of the center of mass.
Regarding the optimization procedure, a natural evolution is the
implementation of a multi-objective optimization which considers at
the same time both mass and deflection at a given point, for example. This can be
done in a Bayesian setting, accounting for high dimensional input
parameter space~\cite{binois2021survey,daulton2021multi}, but also
other approaches should be considered, such as genetic algorithms
enhanced by active subspaces~\cite{demo2020asga, demo2021hull}.

\section*{Acknowledgements}
This work was partially supported by an industrial Ph.D. grant sponsored by
Fincantieri S.p.A. (IRONTH Project), the project SHip OPtimization
with Reduced Order Methods (SH.OP. ROMs) carried out in the context of
the IRISS initiative by SMACT Competence Center, and partially funded by European
Union Funding for Research and Innovation --- Horizon 2020 Program --- in the
framework of European Research Council Executive Agency: H2020 ERC CoG 2015
AROMA-CFD project 681447 ``Advanced Reduced Order Methods with Applications in
Computational Fluid Dynamics'' P.I. Professor Gianluigi Rozza.

\bibliographystyle{abbrvurl}

\end{document}